\documentclass[10pt,namelimits,sumlimits]{ip-journal}
\usepackage{amssymb,amsmath}
\usepackage[mathscr]{eucal}
\usepackage{color}
\usepackage{enumitem}
\usepackage[utf8]{inputenc}
\usepackage[leftcaption]{sidecap}
\usepackage{tikz}
 \usepackage{tikz-cd}
\usetikzlibrary{matrix,arrows,decorations.pathmorphing}
\usetikzlibrary{positioning}
\usetikzlibrary{matrix,arrows,decorations.pathmorphing, shapes.geometric}

\usepackage{psfrag}

\usepackage{epstopdf}
\usepackage{graphicx}
\usepackage{verbatim}
\usepackage{amsmath,amsthm}
\usepackage{graphics}
\usepackage{color}
\usepackage{epsfig}
\usepackage{fullpage}
\usepackage{amssymb,amsmath}
\usepackage[mathscr]{eucal}

\newtheorem{theorem}[equation]{Theorem}
\newtheorem{lemma}[equation]{Lemma}
\newtheorem{prop}[equation]{Proposition}

\newtheorem{definition}[equation]{Definition}

\theoremstyle{remark}
\newtheorem{remark}[equation]{Remark}

\numberwithin{equation}{section}

\newcommand{\dbold}{{\mathbf{d}}}
\newcommand{\Sph}{\mathbb{S}}

\newcommand{\R}{\mathbb{R}}

\newcommand{\HH}{\mathbb{H}}
\newcommand{\ddiv}{\mathrm{div}}

\newcommand{\trace}{\text{Tr}}

\newcommand{\Runder}{\underline{R}}
\newcommand{\rbar}{\overline{r}}

\newcommand{\nablabar}{\nabla}
\newcommand{\nablasigma}{\nabla^{\Sigma}}

\newcommand{\Phip}{\Phi_p}
\newcommand{\Psip}{\Psi_p}
\newcommand{\Phitilde}{\widetilde{\Phi}}
\newcommand{\varphitilde}{\widetilde{\varphi}}

\begin{document}

\title[]{Sharp Area Bounds for Free Boundary Minimal Surfaces in Conformally Euclidean Balls}
\author[B.~Freidin]{Brian~Freidin}
\author[P.~McGrath]{Peter~McGrath}
\date{}
\address{Department of Mathematics, Brown University, Providence,
RI 02912} \email{bfreidin@math.brown.edu}
\address{Department of Mathematics, University of Pennsylvania, Philadelphia PA 19104} 
\email{pjmcgrat@sas.upenn.edu}
\maketitle
\begin{abstract}
We prove that the area of a free boundary minimal surface $\Sigma^2 \subset B^n$, where $B^n$ is a geodesic ball contained in a round hemisphere $\Sph^n_+$, is at least as big as that of a geodesic disk with the same radius as $B^n$; equality is attained only if $\Sigma$ coincides with such a disk.  More generally, we prove analogous results for a class of conformally euclidean ambient spaces.  This follows work of Brendle and Fraser-Schoen in the euclidean setting. 
\end{abstract}

\section{Introduction}
\label{S:intro}
A  properly immersed submanifold $\Sigma^k\subset \Omega^n$ in a domain of a Riemannian manifold is a \emph{free boundary minimal submanifold} if $\Sigma$ is minimal, $\partial \Sigma \subset \partial\Omega$, and $\Sigma$ intersects $\partial \Omega$ orthogonally.  Such submanifolds are volume-critical among all deformations which preserve the condition $\partial \Sigma \subset \partial\Omega$.

There are a variety of interesting uniqueness theorems for free boundary minimal surfaces.  In the 1980s, Nitsche \cite{Nitsche} showed using Hopf differentials that the only free boundary topological disks in the unit ball $B^3\subset \R^3$ are flat. Souam \cite{Souam} extended this result to encompass free boundary 2-disks in 3-balls in space forms, and Fraser-Schoen \cite{FSII} recently extended this further to encompass free boundary 2-disks in balls  of arbitrary dimension in space forms.   

Another direction of recent interest is to prove sharp area bounds for free boundary minimal surfaces.  Fraser-Schoen proved \cite[Theorem 5.4]{FS1} any  free boundary $\Sigma^2 \subset B^n$ has area at least $\pi$; equality holds precisely when $\Sigma$ is congruent to a disk.    Brendle \cite{Brendle:area} proved more generally that free boundary $\Sigma^k \subset B^n$ of arbitrary dimension satisfy the analogous sharp bound $|\Sigma^k |\geq |B^k|$.  A natural question is whether analogous results hold in ambient spaces of constant but nonzero curvature.

\begin{theorem}
\label{Tsphere}
Let $\Sigma^2 \subset B^n$ be a free boundary minimal surface, where $B^n$ is a geodesic ball contained in a hemisphere of $\Sph^n$.  Then $|\Sigma | \geq  |B^2|$, where $|B^2|$ is the area of a geodesic disk with the same radius as $B^n$.  If equality holds, then $\Sigma$ coincides with some such disk. 
\end{theorem}

In Theorem \ref{Tgen} below, we extend Theorem \ref{Tsphere} to a broader class of ambient manifolds.  We consider spaces $(M,g)$ where $M = \Sph^{n-1} \times [0, \rbar)$, $g = dr^2 + h^2(r) g_{\Sph}$, $g_\Sph$ is the round metric on $\Sph^{n-1}$, and the warping function $h$ satisfies
\begin{enumerate}
\item[(C1)] $h(r) = r v(r^2)$, where $v$ is smooth, positive, and satisfies $v(0) = 1$. \label{rho smooth}
\item[(C2)] $(M,g)$ has positive radial  sectional curvature $K(r)  = - h''(r)/h(r)$. \label{rho inc}
\end{enumerate} 
Let $B_R\subset M$ be the geodesic ball centered at $r = 0$ with radius $R$. 
\begin{theorem}
\label{Tgen}
Suppose $h$ satisfies (C1) and (C2).  There exists $\Runder \in (0, \rbar)$ such that for all $R\in (0, \Runder]$ and all free boundary minimal surfaces $\Sigma^2 \subset B^n_R$,  $|\Sigma| \geq |B^2_R|$, where $|B^2_R|$ is the area of a geodesic disk of radius $R$ and center $r=0$.  If equality holds, $\Sigma$ coincides with some such disk.
\end{theorem}

Theorem \ref{Tgen} applies in particular to balls in $(\R^{3}, e^{-|x|^2/4} \delta)$ -- the space in which self-shrinkers for mean curvature flow are minimal surfaces -- and also to $\Sph^n$, which furnishes another proof (see \ref{Rsphere}) of Theorem \ref{Tsphere}.  We note that the condition (C1) implies that $(M,g)$ is smooth and nondegenerate at $r=0$ (cf. \cite[Theorem 1.4]{Brendle:warped}).

The proofs of Theorems \ref{Tsphere} and \ref{Tgen} are motivated by Brendle's ingenious approach in \cite{Brendle:area}.
There Brendle applies the divergence theorem to a vector field $W$ with the following properties:
\begin{enumerate}[label=(\roman*).]
\item  $W$ is defined on $B^n\setminus \{ y\}$, and has a prescribed singularity at $y\in \partial B^n$.
\item  $W$ is tangent to $\partial B^n$ along $\partial B^n\setminus \{y\}$.
\item $\ddiv_\Sigma W\leq 1$ for any submanifold $\Sigma^k \subset B^n$. 
\end{enumerate}

In the euclidean setting of \cite{Brendle:area}, $W$ is a sum of a radial field with divergence bounded above by $1$ centered at $0$ and a singular field with nonpositive divergence centered at $y$.  The analogous field in the setting of Theorem \ref{Tsphere} unfortunately no longer satisfies (iii).  It turns out however that a judiciously chosen convex combination of fields -- each of which has divergence bounded above by 1 -- can be arranged which satisfies (i)-(iii).

An additional challenge in the context of Theorem  \ref{Tgen}  is that the ambient spaces under consideration are generally not homogeneous.  Nonetheless the definition of $W$ there is motivated by the structure -- when written in stereographic coordinates -- of its counterpart from the proof of Theorem \ref{Tsphere}.  Once the definition is made, (i) and (ii) hold automatically, and (iii) holds when the ambient space satisfies the condition (C2).

The calibration vector field strategy in the sprit of \cite{Brendle:area} appears to be quite flexible and has been used recently by Brendle-Hung \cite{Brendle:area2} to prove a sharp lower bound for the area of a minimal submanifold $\Sigma^k \subset B^n$ passing through a prescribed point $y\in B^n$. 

The approach here is also closely related to work of Choe \cite{ChoeCrelle} and Choe-Gulliver \cite{ChoeGulliverM, ChoeGulliverI} on isoperimetric inequalities for domains on minimal surfaces.  Interestingly, while in that setting (cf. also \cite{Bray}) the geometric inequalities are favorable in a negative curvature background, in the present context positive ambient curvature is essential (see \ref{LPhidiv}) to the proof of Theorems \ref{Tsphere}-\ref{Tgen}. 

In Section \ref{S:sphere}, we prove Theorems \ref{Tsphere} and \ref{Tgen}.  In Section \ref{S:LiYau}, we consider notions of conformal volume as introduced by Li-Yau \cite{Li-Yau} and applied to the free boundary setting by Fraser-Schoen \cite{FS1} to give another proof of Theorem \ref{Tsphere}.
\section{Proof of Theorems \ref{Tsphere}-\ref{Tgen}}
\label{S:sphere}

\subsection*{Warped products and radial vector fields}
\label{ss:warped}

Consider an ambient manifold $(M,g)$, where
\begin{align}
\label{Ewarped}
M= \Sph^{n-1}\times [ 0, \rbar), \qquad
g = dr \otimes dr + h(r)^2 g_{\Sph},
\end{align}
$g_{\Sph}$ is the round metric on the unit $\Sph^{n-1}$, and $h$ satisfies (C1) and (C2).   Note that the setting of \ref{Tsphere} is a special case of that of \ref{Tgen}: given $p\in \Sph^n$, the punctured round unit sphere $\Sph^n\setminus \{-p\}$ is isometric to $M$ as in \eqref{Ewarped}, where $\rbar = \pi$ and $h(r) = \sin r$.

Let $B_R$ be the geodesic ball centered at $r = 0$ with radius $R$, and $|B^2_R|$ be the area of any geodesic disk with center at $r=0$ and radius $R$.  More generally, given $p\in M$, we write $\dbold_p$ for the geodesic distance function from $p$ and define a closed geodesic ball about $p$ by 
\begin{align*}
B_\delta(p) := \left\{ q\in M : \dbold_p(q) \leq \delta \right\}.
\end{align*}
 

Throughout, let $\Sigma^2 \subset B_R$ be a minimal surface. Let $\nabla$ be the covariant derivative on $M$ and $\nablasigma, \ddiv_{\Sigma}$, and $\Delta_{\Sigma}$ respectively be the covariant derivative, divergence, and Laplacian operators on $\Sigma$.  It is convenient  to define $\nabla^\Sigma  r^\perp := \nabla  r- \nabla^\Sigma r$; note that $ \left| \nabla^\Sigma r^\perp\right|^2=1-\left|\nabla^\Sigma r\right|^2$.

Theorems  \ref{Tsphere}-\ref{Tgen} follow from the following  general  argument  which shifts the  difficulty of  the  problem to the  construction of a vector field  with certain  properties.  

\begin{prop}
\label{Pproof}
In the setting of either \ref{Tsphere} or \ref{Tgen}, suppose for each $y\in  \partial B_R$, there exists a vector field $W$ on $B^n_R\setminus \{y\}$ with the following  properties:
\begin{enumerate}[label=\emph{(\roman*).}]
\item As $\dbold_y  \searrow0$,  $W  = - | B^2_R| /(\pi \dbold_y) \nabla  \dbold_y  +o \left( \dbold_y^{-1}\right)$.
\item $W$  is tangent to $\partial B_R$ along $\partial B_R\setminus \{ y\} $.
\item  $\ddiv_\Sigma W\leq  1$  for any minimal surface $\Sigma^2 \subset B_R$, with  equality only  if $\nabla^\Sigma r  = \nabla  r$  on $\Sigma$.
\end{enumerate}
Then the conclusion of  \ref{Tsphere}-\ref{Tgen}  holds. 
\end{prop}
\begin{proof}
Fix $y\in \partial \Sigma$ and $W$ as above.  From the  divergence theorem, the minimality  of $\Sigma$, and (iii), 
\begin{align*}
\left| \Sigma \setminus B_{\epsilon}(y) \right| \geq \int_{\Sigma \setminus B_{\epsilon}(y)} \ddiv_\Sigma W 
= \int_{\partial \Sigma \setminus B_{\epsilon}(y)} \langle W, \eta \rangle + \int_{\Sigma \cap \partial B_{\epsilon}(y)} \langle W, \eta \rangle.
\end{align*}
By the free boundary condition, $\eta = \nabla  r$  on $\partial \Sigma$; using (ii) and letting $\epsilon \searrow 0$, we find
\begin{align*}
|\Sigma| \geq \lim_{\epsilon \searrow 0} \int_{\Sigma \cap \partial B_{\epsilon}(y)} \langle W, \eta \rangle.
\end{align*}
On  $\Sigma \cap  \partial B_\epsilon(y)$,  the free boundary condition implies $\eta  = - \nabla\dbold_y +o(1)$; in  combination with  (i) this implies $\langle W,  \eta \rangle =  |B^2_R|/(\pi \epsilon)  + o\left( \epsilon^{-1}\right)$ on  $\Sigma \cap \partial B_\epsilon(y)$.  The free boundary condition also implies $|\Sigma \cap \partial B_\epsilon(y)| =  \pi \epsilon+o(\epsilon)$; after letting $\epsilon\searrow 0$, we conclude from the  preceding  that  $|\Sigma| \geq  |B^2_R|$.

In the case of equality, (iii) implies that the integral curves in $\Sigma$  of $\nabla^\Sigma  r$ are also integral curves in $M$ of $\nabla r$, namely they are parts of geodesics  passing through  $r=0$.  It follows that  $\Sigma$ is a geodesic disk of radius $R$  passing  through $r=0$.
\end{proof}

\begin{remark}
\label{rhexp} The radial sectional curvature of $M$ is $K(r) = - h''(r)/h(r)$.  From this and (C1), it follows that for small $r\geq 0$,
\begin{align*}
h(r) = r - \frac{K(0)}{3!}r^3 + O(r^4).
\end{align*}
\end{remark}

\begin{definition}
\label{dI}
Define functions $I$ and $\varphi$ on $M$ and a vector field $\Phi$ on $M$ by
\begin{align*}
I(r) = \int_0^r h(s)\, d t, 
\qquad 
\varphi(r) = \frac{ I(r)}{h(r)},
\qquad
\Phi(r) = \varphi(r) \frac{\partial }{ \partial r}, 
\end{align*}
where we define $\varphi(0) = 0$. 
\end{definition}
\begin{remark}
\label{rphi}
By \ref{rhexp}, $\varphi$ and $\Phi$ are clearly $C^1$. 
\end{remark}

\begin{lemma}
\label{LPhidiv}
$\displaystyle{ \ddiv_\Sigma \Phi = 1+ 2 h^{-2}\left( \int_0^r I h'' \, dt\right) \left| \nabla^\Sigma r^\perp\right|^2}.$
\end{lemma}
\begin{proof}
Clearly $\nabla I = h\,  \partial_r$, and as in \cite[Lemma 2.2]{Brendle:warped}, $\nabla^2 I =h'(r) g$.
 Consequently, we compute 
 \begin{align}
\label{Ehess}
 \nabla^2 r = \frac{h'}{h} \left( g - dr \otimes dr\right).
 \end{align}
 It follows that $\Delta_\Sigma r =( h'/h) (2 - \left| \nabla^\Sigma r \right|^2)$. 
Then
\begin{align*}
\ddiv_{\Sigma}\Phi &= \varphi' \left| \nabla^\Sigma r\right|^2 + \varphi \Delta_\Sigma r\\
&= \varphi'\left| \nabla^\Sigma r\right|^2 + \varphi \frac {h'}{h} ( 2 - \left| \nabla^\Sigma r\right|^2)\\
&= \varphi'+ \varphi \frac{h'}{h} + \left( \varphi \frac{h'}{h} - \varphi'\right) \left|\nabla^\Sigma r^\perp\right|^2\\
&= 1+ h^{-2}\left({2 I h' - h^2}\right) \left|\nabla^\Sigma r^\perp\right|^2\\
 &=1+ h^{-2} \left( 2\int_{0}^r (I h')' \, dt - h^2\right) \left|\nabla^\Sigma r^\perp\right|^2\\
&=1+2 h^{-2}\left( \int_0^r I h'' \,dt\right) \left|\nabla^\Sigma r^\perp\right|^2.
\end{align*}
\end{proof}

\subsection*{The round sphere $\Sph^n$}  Given $p\in \Sph^n$, recall that $\Sph^n\setminus \{-p\}$ is isometric to $M$ as in \eqref{Ewarped},  where $\rbar = \pi$ and $h(r)=\sin r$.  Let $I$ and $\varphi$ be as in \ref{dI} and compute
\begin{align}
\label{Ephipsi}
I(r)  = 1-  \cos r, \qquad
 \varphi(t) = \frac{1-\cos t}{\sin t} =  \tan \frac{t}{2}.
 \end{align}

\begin{definition}
\label{Dpsi}
Given $p\in \Sph^n$, define vector fields $\Phi_p$ and $\Psi_p$ on respectively $\Sph^n\setminus \{-p\}$ and $\Sph^n \setminus \{ p\}$ by 
\begin{align*}
\Phi_p =\left( \varphi \circ \dbold_p\right) \nablabar \dbold_p, \qquad \Psi_p = \Phi_{-p}.
\end{align*}
\end{definition}

\begin{lemma}
\label{Lphidiv}
Given $p\in \Sph^n$, the following hold.
\begin{enumerate}[label=\emph{(\roman*).}]
\item $\ddiv_\Sigma \Phip \leq 1$ and $\ddiv_\Sigma \Psip\leq 1$.
\item $\Psip =\left( \psi \circ \dbold_p\right) \nablabar \dbold_p$, where $ \psi(t) :=  -\cot(t/2)$.
\end{enumerate}
\end{lemma}
\begin{proof}
(i) follows immediately from \ref{LPhidiv}, \eqref{Ephipsi}, and \ref{Dpsi}.  For (ii), denote $r = \dbold_p$ and compute
\begin{align*}
\Psi_p  &= \tan\left( (\pi -r)/2\right) \nabla  \dbold_{-p}   = -\cot(r/2) \nabla r.
\end{align*}
\end{proof}
Now fix $p\in \Sph^n$, $R\leq \pi/2$, and $y \in \partial B_R$, where $B_R := B_R(p)$. 

\begin{definition}
\label{dW}
Define a vector field on $B_R \setminus \{ y\}$ by 
\begin{align*}
W = \left( \cos R\right)   \Phi_p + (1-\cos R)  \Psi_y.
\end{align*}
\end{definition}

\begin{lemma}
\label{lprods}
With the hypotheses of Theorem \ref{Tsphere}, $W$ satisfies \ref{Pproof}.(i)-(iii).
\end{lemma}
\begin{proof}
Note that $|B^2_R| = 2\pi I(R)   = 2\pi (1-\cos R)$.  (i) then follows from \ref{dW}  by expanding $\Psi_y$ using \ref{Lphidiv}.(ii).  For (ii), compute for $x\in \partial B_R$ 
\begin{align*}
\langle W, \nablabar \dbold_p \rangle  &=\left( \cos R\right) \langle  \Phi_p, \nabla \dbold_p \rangle 
+ (1-\cos R) \langle  \Psi_y, \nabla \dbold_p \rangle\\
&= \cos R \tan\frac{R}{2}  - (1 - \cos R) \cot\frac{\dbold_y}{2} \cos \theta,
\end{align*}
where $\theta$ is the angle in the geodesic triangle $pxy$ subtended at $x$.  By the spherical law of cosines, 
\begin{align*}
\cos \theta = \frac{1-\cos \dbold_y}{\sin \dbold_y} \cot R = \tan \frac{\dbold_y}{2} \cot R.
\end{align*}
It follows that along $\partial B_R$,
\begin{align*}
\langle W, \nabla \dbold_p \rangle =  \cos R \tan \frac{R}{2} -  (1-\cos  R)  \cot  R =  0.
\end{align*} 
(iii) follows immediately from \ref{Lphidiv}.(i), that $R\leq \pi/2$, and \ref{LPhidiv}.  This proves theorem \ref{Tsphere}.
\end{proof}

\subsection*{Proof of Theorem \ref{Tgen}}

It is convenient to use coordinates in which $(M, g)$ is conformal to a euclidean ball $( B^n, \rho^2 \delta)$:

\begin{definition}[Conformal coordinates]
\label{Dconf}
Define a parameter $s=s(r)$ and a function $\rho$ by requesting that  $s(0)=0$ and that
\begin{align}
\label{Econf}
\frac{ds}{dr} = \frac{s(r)}{h(r)}, \qquad h(r) = s\rho(s).
\end{align}
\end{definition}

In these coordinates, $g = \rho^2(s) \delta$, where $\delta$ is the euclidean metric on $\R^n$.  Throughout, we write $\langle\, , \rangle_\delta$ for the metric $\delta$, $|\cdot |_\delta$ for the induced norm, and $s: = |x|_\delta$.  When there is no room for confusion, we abbreviate $\rho = \rho(s)$.  As before we fix $y\in \partial B_R$.

\begin{definition}
\label{dV}
Define a vector field $V$ on $B_R \setminus \{y\}$ by requesting that in conformal coordinates
\begin{align*}
V = \frac{1}{\rho^2(s)} \frac{x-y}{|x-y|_{\delta}^2}.
\end{align*}
\end{definition}
To compute $\ddiv_\Sigma V$ we need the following. 

\begin{lemma}
\label{Ldiv}
Suppose $Y$ is a smooth vector field on a domain of M. Then
	\begin{align*}
	\ddiv_{\Sigma}\left(Y/\rho^2\right) = \frac{1}{\rho^2} \ddiv_{\Sigma, \delta} Y+ 
	\frac{2 \rho'}{ s\rho^{3}} \langle Y, x^\perp\rangle_\delta,
	\end{align*}
	where $\ddiv_{\Sigma, \delta}$ is the divergence operator on $\Sigma$ computed with respect to $\delta$, $x=h(r)\nabla r$ is the position vector field in conformal coordinates, and $x^\perp$ is the part of $x$ orthogonal (with respect to  either $\delta$ or $g$) to $\Sigma$ at $x$.
\end{lemma}
\begin{proof}
A straightforward calculation using the Koszul formula shows that
	\begin{align*}
	\ddiv_{\Sigma} Y = \ddiv_{\Sigma, \delta} Y + \frac{2}{\rho} Y(\rho).
	\end{align*}
Replacing $Y$ by $ Y/\rho^2$ in the above, we find (where below $\nabla^{\Sigma, \delta}$ is the  euclidean connection on $\Sigma$)
\begin{align*}
\ddiv_{\Sigma, g} \left(Y/\rho^2\right) &= \frac{1}{\rho^2} \ddiv_{\Sigma, \delta} Y + \left \langle Y, \nabla^{\Sigma, \delta} \frac{1}{\rho^2} \right \rangle_\delta + \frac{2}{\rho^{3}} Y(\rho)\\
&= \frac{1}{\rho^2}\ddiv_{\Sigma, \delta} Y - \frac{2}{\rho^{3}} \left \langle Y , \nabla^{\Sigma, \delta} \rho\right \rangle_\delta + \frac{2}{\rho^{3}} Y(\rho)\\
&=  \frac{1}{\rho^2} \ddiv_{\Sigma, \delta}Y + \frac{2 \rho'}{ s\rho^{3}} \langle Y, x^\perp\rangle_\delta.
\end{align*}\end{proof}

\begin{lemma} 
\label{Ldiv2}
\begin{enumerate}[label=\emph{(\roman*).}]
\item
$\displaystyle{
\ddiv_\Sigma V \geq
	 - \frac{1}{2}\left( \frac{h' -1}{h}\right)^2\left|\nabla^\Sigma r^\perp\right|^2.}$
	 
\item On $\partial B_R$, $\langle V, \nabla r \rangle = 1/(2 h(R)) $.
\end{enumerate}
\end{lemma}
\begin{proof}
 Using \ref{Ldiv}, we compute in conformal coordinates:
\begin{align*}
 \ddiv_{\Sigma} V &= 
 \frac{1}{\rho^2} \ddiv_{\Sigma, \delta}\left( \frac{x-y}{|x-y|_\delta^2}\right) + \frac{2 \rho'}{s \rho^3} \left \langle \frac{x-y}{|x-y|_\delta^2}, x^\perp\right\rangle_\delta\\
 &= \frac{2}{\rho^2} \frac{|(x-y)^\perp|_\delta^2}{|x-y|_\delta^4} + \frac{2 \rho'}{s \rho^3} \left \langle \frac{x-y}{|x-y|_\delta^2}, x^\perp\right\rangle_\delta\\
 &= 2  \left| \frac{1}{\rho|x-y|_\delta^2} y^\perp - \left( \frac{1}{\rho |x-y|_\delta^2} + \frac{\rho'}{2s \rho^2} \right) x^\perp \right|_\delta^2 - \frac{1}{2}\left(\frac{\rho'}{s\rho^2}\right)^2 |x^\perp|_\delta^2.
 \end{align*}
 Using \eqref{Econf}, note that $\nabla r  = x/ ( s\rho(s))$ and compute that  $|x^\perp|^2_\delta = s^2 \left|\nabla^\Sigma r^\perp\right|^2$.  Differentiating $h(r)  = s\rho(s)$,  we see that $\rho'(s)/\rho(s)= (h'-1)/s$ and (i) follows. 
  For (ii), note that $2\langle x - y, x\rangle_\delta = |x-y|_\delta^2$ when $|x|_\delta = |y|_\delta$.
  Using this and \eqref{Econf} we have
 \begin{align*}
 \langle V, \nabla r \rangle =\left \langle \frac{x-y}{|x-y|^2_\delta}, \frac{x}{s \rho(s)}\right\rangle_\delta
 = \frac{1}{2 h(R)}.
 \end{align*}

\end{proof}

Define now a vector field $W$ on $B^n_R\setminus\{ y\}$ by 
\begin{align}
\label{EW}
W=  \Phi - 2 I(R) V.
\end{align}

\begin{lemma}
\label{LWtan}
With the hypotheses of Theorem \ref{Tgen},  $W$ satisfies (i)-(iii) of \ref{Pproof}.
\end{lemma}
\begin{proof}
Recall  that $|B^2_R|=2\pi I(R)$.  (i) then follows from \eqref{EW}, \ref{dV}, and \eqref{Econf} by taking $x\rightarrow y$.
For (ii), when $r=R$, we have by \ref{dI} and \ref{Ldiv2}.(ii),
\begin{align*}
\left \langle W, \nabla r\right \rangle = \frac{I(R)}{h(R)}- 2I(R)  \langle V, \nabla r\rangle = 0.
\end{align*}
For (iii), we have by \ref{LPhidiv} and \ref{Ldiv2}.(i),
\begin{equation}
\label{Eerror}
\begin{aligned}
\ddiv_\Sigma W &\leq 1 +\left(2\int_0^r I h'' dt + I(R)(h' -1)^2\right) h^{-2} \left|\nabla^\Sigma r^\perp\right|^2\\
&:= 1+ (*)h^{-2}\left|\nabla^\Sigma r^\perp\right|^2.
\end{aligned}
\end{equation}
By straightforward asymptotic expansions using \ref{rhexp} and \ref{dI}, it follows that 
\begin{align*}
(*) = \left( - 1 + K(0) I(R)\right) \frac{K(0)}{4} r^4 + O(r^5).
\end{align*}
Since $I$ is increasing, $K(0)>0$ (recall (C2)), and $I(0) = 0$, (iii) follows by choosing $\Runder$ small enough.
\end{proof}

\begin{remark}
\label{Rsphere}
The conclusion of Theorem \ref{Tgen} holds in particular in the following situations:
\begin{itemize}
\item When $h(r) = \sin r$ and $\rbar = \pi$, i.e. $(M, g)$ is isometric to a  punctured round sphere.  A lengthy computation using stereographic  coordinates shows that $W$ as in \eqref{EW} coincides with the field $W$ in the proof of Theorem \ref{Tgen}.  In particular, one may take $\Runder = 1$, and in that case $B_{\Runder}$ corresponds via stereographic projection to a closed hemisphere of $\Sph^n$. 
\item When $h(r)= r$ and $\rbar =\infty$, i.e. $(M, g)$ is isometric to $(\R^n, \delta)$; then the error term $(*)$ in \eqref{Eerror} is identically $0$ and $W$ coincides with the vector field $W$ defined by Brendle in \cite{Brendle:area} (when $k=2$). 
\item When $M=( \R^3, g:= e^{-|x|^2/4} \delta)$; by a change of coordinates as in \eqref{Econf} $M$ can be brought into the form of \eqref{Ewarped}.  Moreover, as in \cite{CM}, the scalar curvature of $M$ satisfies
\begin{align*}
\text{Scal}_g = e^{\frac{|x|^2}{4}}\left( 3 - \frac{|x|^2}{8}\right),
\end{align*}
which is positive when $|x|^2< 24$.  Direct calculation of $\ddiv_{\Sigma} W$ as in the proof of \ref{Ldiv2} (for which we omit the details) shows that $\Runder$ can be taken as large as the positive root of $8+ r^4 - e^{r^2/2}(8-4r^2+r^4)$, or approximately $\Runder \approx 1.546$. 
\end{itemize}
\end{remark}

\section{Conformal Volume}
\label{S:LiYau}
Li-Yau proved \cite{Li-Yau} that closed minimal surfaces $\Sigma^2 \subset \Sph^n$ maximize their area in their conformal orbit.  Using similar ideas, Fraser-Schoen \cite[Theorem 5.3]{FS1} showed that the boundary $\partial \Sigma$ of a free boundary minimal surface $\Sigma^2 \subset B^n$ maximizes its length in its conformal orbit.  From this they deduced that every free boundary minimal surface $\Sigma \subset B^n$ satisfies $|\Sigma | \geq \pi$.  A modification of their argument leads to another proof of Theorem \ref{Tsphere}.

\begin{proof}[Second proof of Theorem \ref{Tsphere}]  We follow the notation of \cite[Theorem 5.3]{FS1}.
Let $\mathring{A} = A - \frac{1}{2} \left( \trace_g A\right) g$ be the trace-free second fundamental form.  It is well known that $\| \mathring{A} \|^2\,  d S_g$ is a pointwise conformal invariant in dimension two.  Next, note that $2 \| \mathring{A} \|^2 = H^2 - 4 \kappa_1 \kappa_2$, where $\kappa_1$ and $\kappa_2$ are the principal curvatures of $\Sigma$.  The Gauss equation implies $K_\Sigma = 1 + \kappa_1\kappa_2$.  It follows that
\begin{align*}
2\| \mathring{A} \|^2  = H^2 - 4 K_\Sigma + 4.
\end{align*}
Hence, for each conformal transformation $F: B^n_R \rightarrow B^n_R$,
\begin{align*}
\int_{\Sigma} \left( H^2 - 4 K_\Sigma + 4 \right) \, dS = \int_{F(\Sigma)}\left( H^2 - 4 K_{F(\Sigma)} + 4 \right) \, d\widetilde{S},
\end{align*}
where $dS$ and $d\widetilde{S}$ denote the surface measures on $\Sigma$ and $F(\Sigma)$.
Using the minimality of $\Sigma$, we conclude
\begin{align*}
\int_{\Sigma} \left( 1 - K_{\Sigma} \right) \, dV \geq \int_{F(\Sigma)} \left( 1- K_{F(\Sigma)}\right) \, d\widetilde{S}.
\end{align*}
By the Gauss-Bonnet theorem,
\begin{align*}
\int_{\Sigma} K_\Sigma \, dS  &= 2\pi \chi(\Sigma) - \int_{\partial \Sigma} \kappa \, d s,\\
\int_{F(\Sigma)}\! \widetilde{K}_{F(\Sigma)} \, d \widetilde{S} &= 2\pi \chi(F(\Sigma)) - \int_{\partial F(\Sigma)} \widetilde{\kappa} \, d \tilde{s}.
\end{align*}
Combining with the above, we find
\begin{align*}
|\Sigma |- 2\pi \chi(\Sigma) + \int_{\partial \Sigma}\kappa \, ds \geq |F(\Sigma)| - 2\pi \chi(F(\Sigma)) + \int_{\partial F(\Sigma)} \widetilde{\kappa} \, d \tilde{s} , 
\end{align*}
hence
\begin{align}
\label{Econf2}
|\Sigma | +  \int_{\partial \Sigma}\kappa\, d s \geq |F(\Sigma)|+ \int_{\partial F(\Sigma)} \widetilde{\kappa}\, d \tilde{s} .
\end{align}

Let $\gamma$ be a local parametrization of $\partial \Sigma$.  Let $\eta$ be the inward pointing conormal vector field along $\partial \Sigma$ and let $D_s$ be the covariant derivative along $\gamma$.  Then
\begin{align*}
\kappa = \left \langle D_s \dot{\gamma}, \eta \right\rangle = - \left \langle \dot{\gamma}, D_s \eta\right \rangle 
= \left \langle \dot{\gamma} , \nabla_{\dot{\gamma}} \nabla r\right\rangle
=( \nabla^2 r) ( \dot{\gamma}, \dot{\gamma}).
\end{align*}
From \eqref{Ehess}, $\nabla^2 r = \cot r ( g - dr \otimes dr )$; consequently, 
$\kappa = \cot R$. 
Thus, \eqref{Econf2} implies
\begin{align*}
\left(\cot R\right) |\partial\Sigma| + |\Sigma| \geq \left(\cot R\right) |\partial F(\Sigma)| + |F(\Sigma)|.
\end{align*}
Taking a sequence of conformal transformations which concentrates away from a point $y \in \partial \Sigma$ (see \cite[Proposition 1.2]{Souam}, \cite[Fact 2]{Li-Yau} or \cite[Remark 5.2]{FS1} for more details) establishes that
\begin{align}
\label{Elastbound}
\left(\cot R\right) |\partial\Sigma|+ |\Sigma| &\geq  \left(\cot R\right) |\partial B^2_R| + |B^2_R|.
\end{align}
Applying the divergence theorem to $\Phi_p$ (recall \ref{Ephipsi}-\ref{Lphidiv}), we have
\begin{align}
\label{Eeasyarea}
|\Sigma | \geq \tan(R/2) | \partial \Sigma|, \qquad
|B^2_R| = \tan(R/2) |\partial B^2_R|. 
\end{align}
The conclusion then follows from combining \eqref{Elastbound} with \eqref{Eeasyarea}, using that $R\leq \pi/2$.
\end{proof}

\nocite{*}
\bibliographystyle{plain}
\bibliography{bibliography}
\end{document}